\pdfoutput=1
\documentclass{amsart}
\usepackage{graphicx,overpic}
\usepackage{mathptmx}      
 \usepackage{amsmath}
 \usepackage{amsfonts}
 \usepackage{amssymb}
 \usepackage{amsthm}

 \newtheorem{theorem}{Theorem}
 
 \newtheorem{lemma}[theorem]{Lemma}
 \newtheorem{corollary}[theorem]{Corollary}
 
 \theoremstyle{definition}
 
 \newtheorem{conjecture}[theorem]{Conjecture}
\numberwithin{theorem}{section}

\newcommand{\N}{\mathbb{N}}
\newcommand{\Z}{\mathbb{Z}}

\newcommand{\R}{\mathbb{R}}

\DeclareMathOperator{\blocks}{blocks}
\DeclareMathOperator{\gaps}{gaps}
\DeclareMathOperator{\Deg}{deg}
\newcommand{\size}[1]{\left| #1 \right|}

\def\bdc{\mathbf{c}}

\def\bdx{\mathbf{x}}

\def\calD{{\mathcal D}}
\def\calF{{\mathcal F}}
\def\calP{{\mathcal P}}

\def\Min{\operatorname{Min}}
\def\sdepth{\operatorname{sdepth}} 
\def\Supp{\operatorname{Supp}}

\def\depth{\operatorname{depth}}

\def\Set#1{\left\{ #1\right \}}
\def\ceil#1{\left\lceil #1 \right\rceil}
\def\floor#1{\left\lfloor #1 \right\rfloor}

\begin{document}
\date{21 October 2009}
\title{On the Stanley Depth of Squarefree Veronese Ideals}
\author{Mitchel T.~Keller}
\address{School of Mathematics\\ Georgia Institute of Technology\\
  Atlanta, GA 30332-0160 U.S.A.}
\email{keller@math.gatech.edu}

\author{Yi-Huang Shen} \address{Department of Mathematics\\ University
  of Science and Technology of China\\ Hefei, Anhui, 230026, China}
\email{yshen.math@gmail.com}
 
\author{Noah Streib}
\address{School of Mathematics\\ Georgia Institute of Technology\\
  Atlanta, GA 30332-0160 U.S.A.}
\email{nstreib3@math.gatech.edu}

\author{Stephen J.~Young}
\address{School of Mathematics\\ Georgia Institute of Technology\\
  Atlanta, GA 30332-0160 U.S.A.}
\curraddr{Department of Mathematics\\ University of California, San
  Diego\\ La Jolla, CA 92093-0112 U.S.A.}
\email{s7young@math.ucsd.edu}




\begin{abstract}
  Let $K$ be a field and $S=K[x_1,\dots,x_n]$.  In 1982, Stanley
  defined what is now called the Stanley depth of an $S$-module $M$,
  denoted $\sdepth(M)$, and conjectured that $\depth(M) \le
  \sdepth(M)$ for all finitely generated $S$-modules $M$.  This
  conjecture remains open for most cases.  However, Herzog, Vladoiu
  and Zheng recently proposed a method of attack in the case when $M =
  I / J$ with $J \subset I$ being monomial $S$-ideals.  Specifically,
  their method associates $M$ with a partially ordered set.  In this
  paper we take advantage of this association by using combinatorial
  tools to analyze squarefree Veronese ideals in $S$.  In particular,
  if $I_{n,d}$ is the squarefree Veronese ideal generated by all
  squarefree monomials of degree $d$, we show that if $1\le d\le n <
  5d+4$, then $\sdepth(I_{n,d})=
  \floor{\binom{n}{d+1}\Big/\binom{n}{d}}+d$, and if $d\geq 1$ and
  $n\ge 5d+4$, then $d+3\le \sdepth(I_{n,d}) \le
  \floor{\binom{n}{d+1}\Big/\binom{n}{d}}+d$.  
\end{abstract}
\maketitle

\section{Introduction}

Let $K$ be a field and $S=K[x_1,\dots,x_n]$ be the polynomial ring in
$n$ variables with coefficients from $K$. The ring $S$ is naturally
equipped with a $\Z^n$-grading. A \emph{Stanley decomposition} of a
finitely generated $\Z^n$-graded $S$-module $M$ is a representation as
a finite direct sum of $K$-vector spaces
\[
\calD : M=\bigoplus_{i=1}^m u_i K[Z_i],
\]
where each $u_i\in M$ is homogeneous, and $Z_i$ is a subset of
$\Set{x_1,\dots,x_n}$. The \emph{Stanley depth} of this decomposition
is $\sdepth(\calD):=\min\Set{|Z_i| \mid 1\le i\le m}$, where $|Z_i|$
is the cardinality of $Z_i$. The \emph{Stanley depth of $M$} is
defined to be
\[
\sdepth(M):=\max\Set{\sdepth(\calD) \mid \text{ $\calD$ is a Stanley decomposition of $M$}}.
\]
In \cite{MR666158}, Stanley conjectured that $\depth(M) \le
\sdepth(M)$ for all finitely generated $\Z^n$-graded $S$-modules $M$.
This conjecture has been confirmed in several special cases, for
instance, when the module $M$ is {clean} in the sense of Dress
\cite{MR1239277}, that is, it allows a homogeneous prime filtration
$\calF$ with $\Supp(\calF)=\Min(M)$. The Stanley-Reisner ring
$K[\Delta]$ of a simplicial complex $\Delta$ is clean if and only if
$\Delta$ is shellable. It was proved in \cite{popescu-2008} that
Stanley's conjecture also holds for $M=S/I$, where $I\subset S$ is a
monomial ideal and $\dim (S) \le 5$. For other recent developments,
see for example
\cite{MR1958008,MR1958009,arXiv:0906.1303,MR2366164,arxiv.0712.2308}.

However, Stanley's conjecture still remains open, partly because of
the difficulty of computing the Stanley depth of $M$.  Herzog, Vladoiu
and Zheng proposed a method of attack in the special case where
$M=I/J$ with $J\subset I$ being monomial $S$-ideals in
\cite{arxiv.0712.2308}. They associate $I/J$ with a poset $P_{I/J}$,
and thus the computation of $\sdepth(I/J)$ is reduced to a
corresponding computation on the poset $P_{I/J}$. We briefly review
their method for squarefree monomial ideals in
Section~\ref{sec:prelim}. Though a direct implementation of this
method is computationally expensive, it has yielded a wealth of
information; see
\cite{BHKTY,MR2433498,KeYo2009,MR2423489,Okazaki2009,arXiv.org:0805.4461}
for its applications.

In this paper, we investigate the \emph{squarefree Veronese ideal}
$I_{n,d}$, which is generated by all squarefree monomials of degree
$d$ in $S=K[x_1,\dots,x_n]$. The ideal $I_{n,d}$ is the
Stanley-Reisner ideal of the $(d-1)$-skeleton of the $n$-simplex,
hence the Stanley-Reisner ring $S/I_{n,d}$ is Cohen-Macaulay of
dimension $n-d$. The ideal $I_{n,d}$ is a polymatroidal ideal
\cite{MR2215212}, and has linear quotients. The Rees algebras of
squarefree Veronese ideals and their $a$-invariants were studied in
\cite{MR2179232}.

Herzog et al.\ showed in \cite{arxiv.0712.2308} that Stanley's
conjecture holds for both $I_{n,d}$ and $S/I_{n,d}$. However, it is
still difficult to determine $\sdepth(I_{n,d})$. The ideal $I_{n,1}$
is the homogeneous maximal ideal $(x_1,\dots,x_n)$ of $S$. Herzog et
al.\ conjectured in \cite{arxiv.0712.2308} that
$\sdepth(I_{n,1})=\ceil{\frac{n}{2}}$. This formula was later
confirmed by Bir\'o, Howard, Keller, Trotter and Young in \cite{BHKTY}
using combinatorial techniques upon which this paper builds.  The
Stanley depth of $I_{n,d}$ in the other situation remains unexplored
in general, although Cimpoea\c{s} considered a similar problem in
\cite{arXiv:0907.1232}. In this paper, we prove an exact formula for
the Stanley depth of $I_{n,d}$ for certain values of $n$ and $d$ and a
bound for others. Our principal result is the following theorem:

\begin{theorem}\label{thm:main}
  Let $K$ be a field, $S=K[x_1,\dots,x_n]$ be a polynomial ring in $n$
  variables over $K$. Furthemore, let $I_{n,d}$ be the squarefree
  Veronese ideal in $S$, generated by all squarefree monomials of
  degree $d$.
  \begin{enumerate}
  \item If $1\le d\le n < 5d+4$, then the Stanley depth $\sdepth(I_{n,d})= \floor{\binom{n}{d+1}\Big/\binom{n}{d}}+d$.
  \item If $d\geq 1$ and $n\ge 5d+4$, then 
      \[
      d+3\le \sdepth(I_{n,d}) \le \floor{\binom{n}{d+1}\Big/\binom{n}{d}}+d.
      \]
  \end{enumerate}
\end{theorem}

We begin the paper by reviewing the method of Herzog et al.\ for
associating a poset with a monomial ideal. We also introduce our
notation and prove some preliminary results in
Section~\ref{sec:prelim}. Section~\ref{sec:proof} is devoted to
proving Theorem~\ref{thm:main} in three stages.

\section{Background, notation, and preliminary results}
\label{sec:prelim}

\subsection{From ideals to interval partitions}

By definition, the ideal $I_{n,d}$ is squarefree. The method of Herzog
et al.\ for determining the Stanley depth of a squarefree monomial
ideal $I$ using posets can be summarized in the following way. Let
$G(I)=\Set{v_1,\dots,v_m}$ be the unique set of minimal monomial
generators of $I$. For each $\bdc=(c(1),\dots,c(n))\in \N^n$, denote
$\bdx^\bdc=\prod_i x_i^{c(i)}$. The monomial $\bdx^\bdc$ is squarefree
precisely when $c(i)\in \Set{0,1}$ for $1\le i\le n$. By definition,
all the $v_i$ are squarefree. If $\bdx^\bdc$ is squarefree, let
$\Supp(\bdx^{\bdc})=\Set{i\mid c(i)=1}\subseteq [n]:=\{1,\dots,n\}$ be
the support of $\bdx^\bdc$. For the squarefree ideal
$I=(v_1,\dots,v_m)$, there is a poset with ground set
\[
P_I=\Set{ C \subseteq [n] \mid  C \text{ contains $\Supp(v_i)$ for some $i$}}
\]
partially ordered by inclusion. The poset $P_I$ is a subposet of the
Boolean algebra of subsets of $[n]$; it consists of all supersets of
the support of the generators of $I$. For every $A, B\in P_I$ with
$A\subseteq B$, define the interval $[A,B]$ to be $\Set{C \in P_I \mid
  A \subseteq C \subseteq B}$. Let $\calP: P_I=\bigcup_{i=1}^r [C_i,
D_i]$ be a partition of $P_I$, and for each $i$, let $\bdc_i \in \N^n$
be the tuples such that $\Supp(\bdx^{\bdc_i})=C_i$. There is then
a Stanley decomposition $\calD(\calP)$ of $I$:
\[
\calD(\calP): I=\bigoplus_{i=1}^r \bdx^{\bdc_i}K[\Set{x_j \mid j\in D_i}].
\]
It is clear that
$\sdepth(\calD(\calP))=\min\Set{|D_1|,\dots,|D_r|}$. Most importantly,
Herzog et al.\ showed in~\cite{arxiv.0712.2308} that
if $I$ is a squarefree monomial ideal, then
\[
\sdepth(I)=\max\Set{\sdepth(\calD(\calP))\mid \text{$\calP$ is a partition of $P_I$}}.
\]

Before introducing our combinatorial approach to this problem, we
establish some bounds on the Stanley depth of squarefree Veronese
ideals. We also prove other results which show that our approach to
proving Theorem~\ref{thm:main} can focus on values of $n$ having a
particular form.

\begin{lemma}\label{lem:min-sdepth}
    If $I$ is a squarefree monomial ideal and $G(I)$ is the minimal monomial generating set of $I$, then $ \sdepth(I)\ge \min\Set{\deg u \mid u\in G(I)}$.
\end{lemma}

\begin{proof}
    This follows directly from the method of Herzog et al.\
\end{proof}

\begin{lemma}
    \label{LEM2}
    For integers $1\leq d\leq n$, $\sdepth(I_{n,d})\le
    \binom{n}{d+1}/\binom{n}{d}+d$.
\end{lemma}

\begin{proof}
  Since $I_{n,d}$ is generated by monomials of degree $d$, by
  Lemma~\ref{lem:min-sdepth} we have $k:=\sdepth(I_{n,d})\ge d$. The
  poset $P_{I_{n,d}}$ has a partition $\calP: P_I=\bigcup_{i=1}^r
  [C_i, D_i]$, satisfying $\sdepth(\calD(\calP))=k$. For each interval
  $[C_i,D_i]$ in $\calP$ with $|C_i|=d$, we have $|D_i|\ge
  k$. Furthermore, there are $|D_i|-|C_i|$ subsets of cardinality
  $d+1$ in this interval. Since these intervals are disjoint, counting
  the number of subsets of cardinality $d+1$, we have the inequality:
    \[
    \binom{n}{d}(k-d)\le \binom{n}{d+1}.
    \]
    After simplification, we get the desired estimate.
\end{proof}

The following corollary is now immediate.

\begin{corollary}
    If $d\ge \ceil{\frac{n}{2}}$, then $\sdepth(I_{n,d})=d$.
\end{corollary}

We expect the following conjecture, which generalizes our
Theorem~\ref{thm:main}, to be true.

\begin{conjecture}
    \label{CONJ-IND}
    For positive integers $1\leq d\leq n$, $\sdepth(I_{n,d})=
    \floor{\binom{n}{d+1}/\binom{n}{d}}+d$.
\end{conjecture}

This formula is already known to be valid for $d=1$, see
\cite{BHKTY}. A further implementation by computer shows that it is
also valid for $1\le d \le n \le 13$.

The structure of the Boolean algebra of subsets of $[n]$ provides the
following lemma which allows us to use known bounds on the Stanley
depth of a squarefree Veronese ideal to establish other bounds.

\begin{lemma}
  \label{plus1}
  Suppose $d,n,a$ are positive integers such that
  $\sdepth(I_{n,d-1})\ge a$ and $\sdepth(I_{n,d})\ge a+1$. Then
  $\sdepth(I_{n+1,d})\ge a+1$.
\end{lemma}

\begin{proof}
  The poset $P_{I_{n,d-1}}$ has a partition $\calP_1$ with
  $\sdepth(\calD(\calP_1))=a$, while $P_{I_{n,d}}$ has a partition
  $\calP_2$ with $\sdepth(\calD(\calP_2))=a+1$. Now
    \[
    \calP: P_{I_{n+1,d}} = \left(\bigcup_{[C,D]\in \calP_1} \left[C\cup\Set{n+1},D \cup\Set{n+1}\right]\right) \cup \calP_2
    \]

  \item is a partition of $P_{I_{n+1,d}}$ with $\sdepth(\calD(\calP))=a+1$.
\end{proof}

Notice that when $\binom{n}{d+1}/\binom{n}{d}=c-1\in \Z$, we have
$n=cd+c-1$. In fact, the following lemma establishes that if we can
prove Conjecture~\ref{CONJ-IND} for such values of $n$, it is true for
all $n$ and $d$, $1\leq d\leq n$.

\begin{lemma} \label{Reduction} Suppose $\sdepth(I_{cd+c-1,d})=d+c-1$
  for all positive integers $c$ and $d$.  Then $\sdepth(I_{n,d}) \geq
  d+c-1$ for all $c$ such that $cd+c-1 \leq n$.
\end{lemma}

\begin{proof}
  The proof proceeds by double induction, first on $d$ and then on
  $n$. For the induction on $d$, notice that in \cite{BHKTY}, Biro et
  al.\ proved that when $d=1$, $\sdepth(I_{n,d}) \geq \lceil n/2
  \rceil$. So if $n \geq cd+c-1 = 2c-1$ then
  \[\sdepth(I_{n,d}) \geq \ceil{ \frac{n}{2}}  \geq \ceil{
    \frac{2c-1}{2}} = c = d+c-1.\]
  So suppose that for all $1 \leq d' \leq d$ and for all $c \geq
  1$, we have for all $n \geq cd' + c -1$ that $\sdepth(I_{n,d})
  \geq d' + c-1$.

  Now consider $I_{n,d+1}$. For this, we use induction on $n$. For the
  base case, suppose $n = d+1$, which is the smallest $n$ can be since
  $n \geq c(d+1)+c-1$ and $c \geq 1$. This case follows from the
  assumption of the lemma with $c=1$. By assumption we have that if $n
  = c(d+1)+c-1$ for some $c$, then $\sdepth(I_{n,d+1}) \geq
  (d+1)+c-1$.  Hence, for all $1 \leq c' \leq c$ for which $n \geq
  c'(d+1)+c'-1$, we have $\sdepth(I_{n,d+1}) \geq c'(d+1)+c'-1$.

  Thus suppose there exists some $c'' \geq 1$ such that
  \[c''(d+1) +c''-1 ~<~ n ~<~ (c''+1)(d+1) + (c''+1) - 1.\] As $n-1
  \geq c''(d+1) + c'' -1$, our inductive hypothesis on $n$ yields
  $\sdepth(I_{n-1,d+1}) \geq (d+1) + c'' -1$. Also, $ n > c''(d+1) +
  c''-1 = c''d + 2c-1$, and therefore we have $n-1 > c''d +c''-1.$
  Thus the inductive hypothesis for $d$ yields $\sdepth(I_{n-1,d})
  \geq d + c'' -1$. Thus, by Lemma \ref{plus1}, $\sdepth(I_{n,d+1})
  \geq (d+1) + c'' -1$, as desired. This concludes both inductions and
  finishes the proof of the lemma.
\end{proof}

Observe that if $c = \floor{\frac{n+1}{d+1}}$ then $cd+c-1 \leq
n$. Thus if $\sdepth(I_{cd+c-1,d}) \geq c+d-1$ for all positive
integers $c$ and $d$, then by Lemma \ref{LEM2} and Lemma
\ref{Reduction} we have that for all $n \geq d$,
\[d + \floor{\frac{n-d}{d+1}} = \floor{\frac{n+1}{d+1}} + d -1 \leq
\sdepth(I_{n,d}) \leq d +\floor{\frac{\binom{n}{d+1}}{\binom{n}{d}}} =
d +\floor{\frac{n-d}{d+1}}.\] Thus, in order to prove Conjecture
\ref{CONJ-IND} it suffices to show that $\sdepth(I_{cd+c-1,d}) \geq
d+c-1$ for all positive integers $c$ and $d$. In other words, we may
restrict our attention to $I_{n,d}$ where $1 \leq d \leq n$ and $n =
cd+c-1$ for some $c \geq 1$.




\subsection{Combinatorial Tools and Definitions}

In this paper, we will frequently refer to the Boolean algebra of
subsets of $[n]$ partially ordered by inclusion as the
\emph{$n$-cube}. A $t$-set is a subset of $[n]$ of size $t$. Recall
that for $A \subseteq B$ in $[n]$, we have defined the {\it interval}
$[A,B]$ as the set of all $C \subseteq [n]$ such that $A \subseteq C
\subseteq B$.  Notice that $[A,B]$ is isomorphic to the
$(|B|-|A|)$-cube. We say that an interval $[A,B]$ {\it covers} a set
$C$ if $C \in [A,B]$. For our purposes, it is helpful to think of a
subset $S$ of $[n]$ by evenly distributing the points $1,2,\dots, n$
around a circle in the plane, with the elements of $S$ depicted as
solid points and the elements of $[n]\backslash S$ depicted as open
points. As such, the circle will have inherent clockwise and
counterclockwise directions.  We will refer to this as the {\it
  circular representation of $[n]$}.

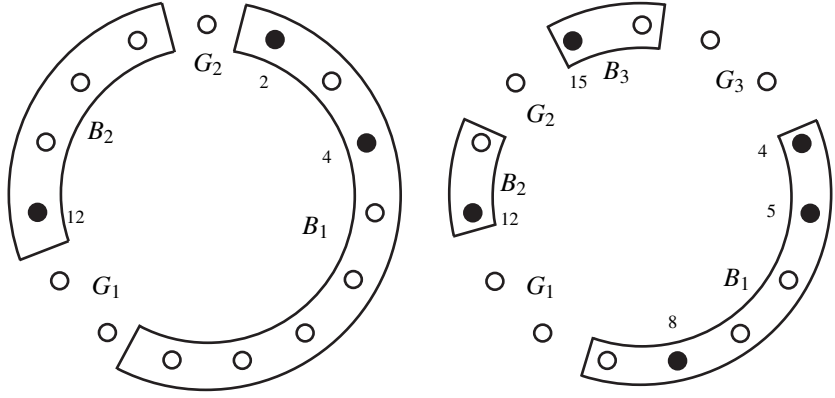
\begin{figure}
\begin{center}
\begin{overpic}{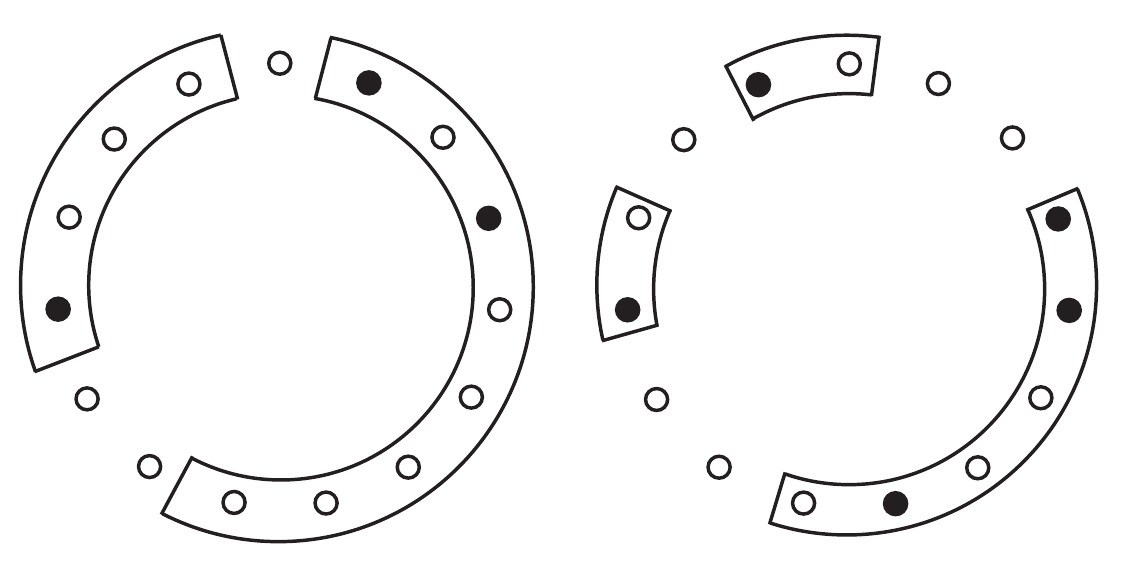}
  \put(31,37){\scriptsize 2}
  \put(38.3,28.3){\scriptsize 4}
  \put(8.5,21.3){\scriptsize 12}
  \put(23.3,39){$G_2$}
  \put(36,20){$B_1$}
  \put(11.5,13){$G_1$}
  \put(11,31){$B_2$}

  \put(89,29){\scriptsize 4}
  \put(90,22){\scriptsize 5}
  \put(78.5,9){\scriptsize 8}
  \put(59,21){\scriptsize 12}
  \put(67,37){\scriptsize 15}
  \put(85,14){$B_1$}
  \put(62,13){$G_1$}
  \put(59,25){$B_2$}
  \put(62,33){$G_2$}
  \put(71,38){$B_3$}
  \put(84,37){$G_3$}
\end{overpic}
\caption{Two examples of block stuctures of $A$ with respect to
  $\delta$. On the left, $A=\{2,4,12\}$ and $\delta=4$. On the right,
  $A=\{4,5,8,12,15\}$ and $\delta=5/2$}
\label{fig:block-structure}
\end{center}
\end{figure}

In order to facilitate our proofs in this paper, we first establish
some definitions. Given the circular representation of $[n]$, a {\it
  block} is a subset of consecutive points on the circle.  For $i,j
\in [n]$ we denote by $[i,j]$ the block starting at $i$ and ending at
$j$ when traversing the circular representation of $[n]$ clockwise.
Let $A \subseteq [n]$ and let $\delta \in \R$ with $\delta \geq 1$ be
a {\it density}.  The {\it block structure of $A$ with respect to
  $\delta$} is a partition of the elements of the circular
representation of $[n]$ into clockwise-consecutive blocks $B_1, G_1,
B_2, G_2, \dots, B_k, G_k$ such that
\begin{itemize}
\item[(i)] the first (going clockwise) element $b_i$ of $B_i$ is in $A$,
\item[(ii)] for all $i \in [k]$, $G_i \cap A = \varnothing$,
\item[(iii)] for all $i \in [k]$, $\delta \cdot |A \cap B_i| - 1 ~<~ |B_i| ~\leq~ \delta \cdot |A \cap B_i|$, 
\item[(iv)] for all $y \in B_i$ such that $[b_i,y] \subsetneq B_i$, $|[b_i,y]|+1 \leq \delta \cdot |[b_i,y] \cap A|.$
\end{itemize}
The purpose of (iii) and (iv) in the definition of block structure
above is to ensure that the density of elements of $A$ in a block is
close to $1/\delta$. For two examples of a block structure, see
Fig.~\ref{fig:block-structure}.

\begin{lemma} For $1 \leq \delta \leq (n-1)/\size{A}$ the block
  structure for a set $A$ on $[n]$ exists and is unique.
\end{lemma}
\begin{proof}
  It is clear for $\delta$ in this range that a block structure can be
  constructed by iteratively expanding blocks.  That is, for each
  element of $A$ assign the next $\floor{\delta}$ elements cyclically
  to the block containing $A$.  If any of these blocks contain another
  element of $A$, merge the two blocks and then add points from
  $[n]\backslash A$ if the number of elements of $A$ in the block
  allows. This process will eventually terminate, yielding a block
  structure, since $\delta \size{A} \leq n-1$.

  Now suppose that $B_1, G_1, \ldots, B_k, G_k$ and $B'_1, G'_1,
  \ldots, B'_{k'}, G'_{k'}$ are two block structures for the set $A$
  with different blocks.  Then there are two blocks, say $B_i$ and
  $B'_j$, which are not the same and overlap.  From (iii) and (iv) it
  is clear that $b_i \neq b'_j$ and so without loss of generality we
  may assume that $b_i$ precedes $b'_j$ in the clockwise ordering.  We
  then note that if $B'_j \setminus B_i$ is nonempty; that is, there
  is an element of $B'_j$ following $B_i$, then
  \begin{align*}
    \size{B_i} &= \size{B_i \setminus B'_j} + \size{B_i \cap B'_j} \\
    &\leq \delta \size{ (B_i \setminus B'_j) \cap A} - 1 + \delta \size{ (B_i \cap
      B'_j) \cap A} - 1 \\
    &= \delta \size{B_i \cap A} - 2,
  \end{align*}
  where the inequality comes from (iv), applied to the initial
  segments of $B_i$ and $B'_j$, respectively.  Note that this
  contradicts (iii), as $\delta\size{B_i \cap A} -1$ is not less than
  $\size{B_i}$.  Furthermore, this implies that $B'_j \subseteq B_i
  \setminus \left\{b_i\right\}$.  But since $b_i \in A$, there must be
  some block $B'_t$ ($i \neq j$) that contains $b_i$.  Clearly, $B_i
  \setminus B'_t$ is nonempty.  Specifically, it contains $B'_j$, and
  so by the above, $B_i \subset B'_t$.  But then $B'_j \subset B_i
  \subset B'_t$, contradicting that the block structure is a
  partition.
\end{proof}

We denote the set $\{B_1, B_2, \dots, B_k\}$ by $\delta$-$\blocks(A)$
and the union $B_1 \cup B_2 \cup \dots \cup B_k$ by
$\mathcal{B}_\delta(A)$.  Each $G_i$ for $i \in [k]$ is called a {\it
  gap.}  We denote the set $\{G_1, G_2, \dots, G_k\}$ by
$\delta$-$\gaps(A)$ and the union $G_1 \cup G_2 \cup \dots \cup G_k$
by $\mathcal{G}_\delta(A)$. Given density $\delta$, let $f_\delta$ be
the function that maps each $A \subseteq [n]$ with $|A|\le
(n-1)/\delta$ to $A\cup \mathcal{G}_\delta(A) \subseteq [n]$.
Throughout this paper we will be concerned with intervals of the form
$[A, f_{\delta}(A)]$.

We will also require the use of some graph theoretic notions. The {\it
  neighborhood} of a vertex $v$ is the set of all vertices adjacent to
$v$ and is denoted $N(v)$. If $S$ is a set of vertices, then $N(S)$
denotes all vertices not in $S$ that are adjacent to a vertex in
$S$. A {\it bipartite graph} is a graph whose vertex set can be
partitioned into sets $V_1$ and $V_2$ such that all edges of $G$ have
one end in $V_1$ and the other in $V_2$.  A {\it matching} in a graph
is a set of edges in which no two edges in the set have a common end.
A {\it complete matching from $V_1$ to $V_2$} in a bipartite graph $G$
with bipartition $(V_1, V_2)$ is a matching containing precisely one
edge incident with each vertex in $V_1$.  A famous result concerning
complete matchings is the following theorem due to Hall.

\begin{theorem}[Hall] \label{thm:hall}
Let $G$ be a bipartite graph with bipartition $(V_1,V_2)$.  Then $G$ has a complete matching from $V_1$ to $V_2$ if and only if $|N(S)| \geq |S|$ for all $S \subseteq V_1$.  
\end{theorem}

In this paper, we will find the following consequence of Hall's
Theorem useful.

\begin{corollary} \label{cor:spec-hall} Let $G$ be a bipartite graph with
  bipartition $V(G) = (V_1,V_2)$. If each vertex in $V_1$ has degree
  $t$ and $t$ is the maximum degree of $G$, then $G$ has a complete
  matching from $V_1$ to $V_2$.
\end{corollary}

The following lemma provides the first step toward proving
Theorem~\ref{thm:main}.

\begin{lemma}\label{lem:c2}
  Let $d\geq 1$ be an integer and $n=2d+1$. Then
  $\sdepth(I_{n,d})=d+1$.
\end{lemma}

\begin{proof}
  It suffices to show that the $d$-sets in $P_{I_{n,d}}$ can be
  covered by intervals $[A,B]$ with $|B|=d+1$. Consider the bipartite
  graph in which $V_1$ consists of the $d$-sets of $[n]$, $V_2$
  consists of the $(d+1)$-sets, and a $d$-set is adjacent to precisely
  the $(d+1)$-sets of which it is a subset. Then every vertex has
  degree $d+1$, so by applying Corollary~\ref{cor:spec-hall} to $V_1$,
  there is a complete matching from $V_1$ to $V_2$. The edges of this
  matching then give rise to intervals of the desired form.  
\end{proof}

Lemma~\ref{lem:c2} is the case of $n=cd+c-1$ with $c=2$. The next
section tackles the more challenging cases of $c=3$ and $c=4$.

\section{Proof of Theorem~\ref{thm:main}}
\label{sec:proof}

To prove Theorem~\ref{thm:main}, it suffices by Lemma~\ref{LEM2},
Lemma~\ref{plus1}, Lemma~\ref{lem:c2} and the proof of
Lemma~\ref{Reduction} to show that if $n=cd+c-1$ and $c\in\{3,4\}$,
then $\sdepth(I_{n,d})=d+c-1$. To do this, we first prove some lemmas
about intervals of the form $[A,f_\delta(A)]$. Our first lemma deals
with the disjointness of these intervals in certain cases.

\begin{lemma} \label{no collisions} Given a positive integer $n$, let
  $A, A' \subseteq [n]$ with $A \not= A'$ and $|A| = |A'|$, and let
  $\delta \in \R$ with $\delta \geq 1$.  If $|f_{\delta}(A)|-|A| \leq
  \delta-1$, then $[A,f_{\delta}(A)]$ does not intersect
  $[A',f_{\delta}(A')]$.
\end{lemma}
\begin{proof}
  Suppose not.  Then $[A,f_{\delta}(A)]$ and $[A',f_{\delta}(A')]$
  have nontrivial intersection.  In particular, $(A \cup A') \in
  [A,f_{\delta}(A)] \cap [A',f_{\delta}(A')]$.  This implies that $a
  \in A' \cup \mathcal{G}_\delta(A')$ for all $a \in A$, and $a' \in A
  \cup \mathcal{G}_\delta(A)$ for all $a' \in A'$.  As $A \not= A'$,
  there exists an $x \in A' \setminus A$ satisfying $x \in
  \mathcal{G}_\delta(A)$.  Let $B_1, G_1, B_2, \dots, B_k, G_k$ be as
  in the definition of the block structure of $A$ and, without loss of
  generality, suppose $x \in G_1$.  Note that
  \[|\mathcal{G}_\delta(A)|~=~|f_{\delta}(A)|-|A|~\leq~\delta-1.\]
  Therefore, the block in $\delta$-$\blocks(A')$ that contains $x$,
  say $B'$, must also contain $b_{2_1}$, the first element of $B_2$.
  As $b_{2_1} \in A$ and $b_{2_1} \notin \mathcal{G}_\delta(A')$ it
  must be that $b_{2_1} \in A'$.  Let $b_{2_2}$ be the next element of
  $B_2\cap A$ found when proceeding clockwise around the circular
  representation of $[n]$.  Clearly $b_{2_2} \in
  B'$, so it must be the case that $b_{2_2} \in
  A'$ as otherwise $b_{2_2} \notin f_{\delta}(A')$.  Proceeding
  clockwise in this manner, we find $A \cap B_2 \subseteq A'$ and
  $B_2 \subseteq B'$.

  Again, using the fact that $|\mathcal{G}_\delta(A)| \leq \delta-1$
  and the fact that $x \in B'$, we find that $b_{3_1}$, the first
  element of $B_3$, is also in $B'$.  Applying the same argument as
  was applied to $B_2$ we find that $A \cap B_2 \subseteq A'$ and $B_3
  \subseteq B'$.

  Finally, proceeding clockwise and using analogous arguments, we
  conclude that $A \cap \mathcal{B}_\delta(A) \subseteq A'$.  However,
  $A \cap \mathcal{B}_\delta(A) = A$, contradicting the fact that $A
  \not\subseteq A'$.
\end{proof}


We also require the following lemma about the size of $f_c(A)$.

\begin{lemma} \label{right size} Let $c,d$ be positive integers and
  let $n = cd+c-1$.  Let $A \subseteq [n]$ be a $d$-set.  Then
  $|f_c(A)| = d+c-1$.
\end{lemma}
\begin{proof}
  Since $c \in \N$, it follows from the definition of the block
  structure of $A$ with respect to the density $c$ that
  $|\mathcal{B}_c(A)| = c|A|$.  Since $A$ is a $d$-set, we find that
  $|\mathcal{G}_c(A)| = n - cd = c-1$.  Thus, $|f_c(A)| = |A| +
  |\mathcal{G}_c(A)| = d+c-1$.  \end{proof}

We now show that the intervals of the form $[A,f_c(A)]$ with $|A|=d$
cover all the $(d+1)$-subsets of $[n]$.

\begin{lemma} \label{d+1} Let $c,d$ be positive integers and let $n =
  cd+c-1$.  Let \[ \mathcal{I} = \{ [A,f_c(A)]\mid A \subseteq
  [n],~|A|=d \}.\] Then each $(d+1)$-subset of $[n]$ is covered by a
  unique element of $\mathcal{I}$.
\end{lemma}
\begin{proof}
Lemma \ref{right size} implies that $|f_c(A)|-|A| = c-1$.  Thus, we may apply Lemma \ref{no collisions}, which implies that each $(d+1)$-set is covered by at most one element of $\mathcal{I}$.  So it suffices to show that the number of $(d+1)$-sets covered by the elements of $\mathcal{I}$ is precisely $\binom{n}{d+1}$.

By Lemma \ref{right size}, each element of $\mathcal{I}$ is a
$(c-1)$-cube.  Thus, each element of $\mathcal{I}$ covers exactly
$c-1$ sets of size $(d+1)$.  Since $|\mathcal{I}| = \binom{n}{d}$, we
find that the elements of $\mathcal{I}$ cover $(c-1)\binom{n}{d}$
sets.  Moreover, since $n=cd+c-1$, it is easy to verify that 
\[	(c-1)\binom{n}{d} = \binom{n}{d+1} \]
as desired.
\end{proof}

The following theorem resolves the first of our two remaining cases,
that of $c=3$.

\begin{theorem} \label{c=3} Let $d$ be a positive integer and
  let $n = 3d + 2$.  Then there exists a partition of the
  subsets of $[n]$ of size at least $d$ into intervals $[A,B]$ such
  that $|B| \geq d+2$.
\end{theorem}
\begin{proof}
  Let $\mathcal{I} = \{[A, f_3(A)]\mid A \subseteq [n], |A| = d\}$.
  By Lemma \ref{right size}, $|f_3(A)| = d+2$ for all $d$-sets $A$.
  By Lemma \ref{d+1}, all $(d+1)$-sets are covered by the elements of
  $\mathcal{I}$.  Since Lemma \ref{right size} implies that
  $|f_3(A)|-|A| = 2$, we can apply Lemma \ref{no collisions} to the
  $d$-sets to show that $I \cap I' = \varnothing$ for all $I \not= I'
  \in \mathcal{I}$.  Moreover, the $(d+2)$-sets that are not covered
  by any element of $\mathcal{I}$ and the subsets of $[n]$ of size at
  least $d+3$ can be covered with trivial intervals.  Therefore we
  have found a partition of the desired type.  \end{proof}

The next lemma provides the final tool required for our proof. It
shows that after selecting the intervals of the form $[A,f_c(A)]$ with
$|A|=d$ for our partition, any uncovered set has all of its supersets
uncovered as well.

\begin{lemma} \label{uncovered top} Let $c,d$ be positive integers
  with $c\geq 2$. Let $n = cd+c-1$ and let $k$ be a positive integer.
  Let $\mathcal{I} = \{[A,f_c(A)]\mid A \subseteq [n], |A|=d\}$.
  Suppose $D_k$ is a $(d+k)$-set that is not covered by any element of
  $\mathcal{I}$.  Then there is no superset of $D_k$ that is covered
  by an element of $\mathcal{I}$.
\end{lemma}
\begin{proof}
  Suppose not. Let $D'$ be a superset of $D_k$ that is covered by some
  $I \in \mathcal{I}$.  By Lemma \ref{right size}, $I$ covers some
  $(d+c-1)$-set, say $D''$, that is also a superset of $D_k$.  Thus,
  it suffices to prove that there is no superset of $D_k$ of size
  $(d+c-1)$ that is covered by an element of $\mathcal{I}$.  To this
  end, let $[D_0, f_c(D_0)] \in \mathcal{I}$ with $D_k \subseteq
  f_c(D_0)$.  Let $X = f_c(D_0) \setminus D_k$.  We call such a
  combination of sets $(X,D_0)$ a {\it pair}.  We call the pair
  $(X,D_0)$ {\it optimal} if, among all pairs, $|X \cap D_0|$ is
  minimized.

  Let $(X^{(0)},D_0^{(0)})$ be an optimal pair. Notice that if $X^{(0)} \cap
  D_0^{(0)} = \varnothing$, then since we know $D_0^{(0)}\subseteq X^{(0)}\cup
  D_k$, we have that $D_0^{(0)} \subseteq D_k \subseteq f_c(D_0^{(0)})$. In this
  case, $[D_0^{(0)}, f_c(D_0^{(0)})]$ covers $D_k$, a contradiction.  Thus
  $|X^{(0)} \cap D_0^{(0)}| \geq 1$.  Consider $x_0 \in X^{(0)} \cap D_0^{(0)}$.
  Let $B \in c$-$\blocks(D_0^{(0)})$ such that $x_0 \in B$. Let $x_1$ be the
  first element in $\mathcal{G}_c(D_0^{(0)})$ counterclockwise from $B$ in
  the circular representation of $[n]$, and let $z_0$ be the last
  element of $B$ (that is, the most clockwise element of $B$). The
  point $z_0$ exists and is distinct from $x$ since for $\delta\geq
  2$, a block cannot end in an element of $A$. Let $x_2,\dots,x_t$ be
  the successive elements of the gaps of $D_0^{(0)}$, indexed
  \emph{counterclockwise} from $x_1$. Fix $s$ as small as possible so
  that $x_{s+1}\in D_k$. Such an $s$ must exist, as otherwise all the
  gap points belong to $X^{(0)}$, so we would have $|X^{(0)}|\geq
  c-1$, but $|X^{(0)}|=d+c-1-(d+k) = c-1-k < c-1$.

  We now define a sequence of pairs $(X^{(i)},D_0^{(i)})$. Let $0\leq
  i\leq s$ and then define $D_0^{(i+1)} = (D_0^{(i)} \setminus
  \{x_i\}) \cup \{x_{i+1}\}$. Notice that $f_c(D_0^{(i+1)}) =
  (f_c(D_0)\setminus\{x_i\}) \cup \{z_i\}$, where $z_i$ is the last
  element of the block of $c$-$\blocks(D_0^{(i)})$ that contains
  $x_i$. We then define $X^{(i+1)} = f_c(D_0^{(i+1)}) \setminus D_k =
  (X^{(i)} \setminus \{x_i\}) \cup \{z_i\}$.  See, for example,
  Fig.~\ref{fig:changing-block-structure}.
  \begin{figure}
    \begin{center}
      \begin{overpic}[scale=0.9]{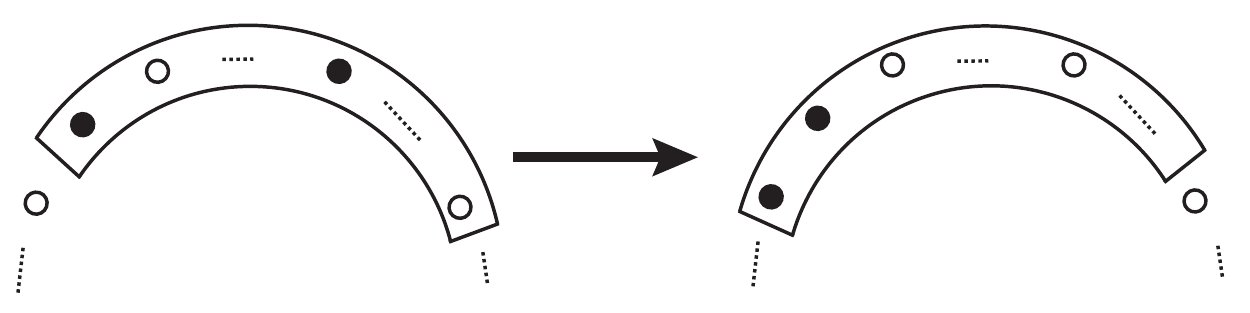}
        \put(5,8){$x_{i+1}$}
        \put(21,24){$B$}
        \put(25,15){$x_i$}
        \put(33,7){$z_i$}
        
        \put(65,7){$x_{i+1}$}
        \put(78,24){$B$}
        \put(84,15){$x_i$}
        \put(92,8){$z_i$}
      \end{overpic}
     \caption{On the left is the block structure of $D_0^{(i)}$, while on the
        right is the block structure of $D_0^{(i+1)}$.}
      \label{fig:changing-block-structure}
    \end{center}
  \end{figure}
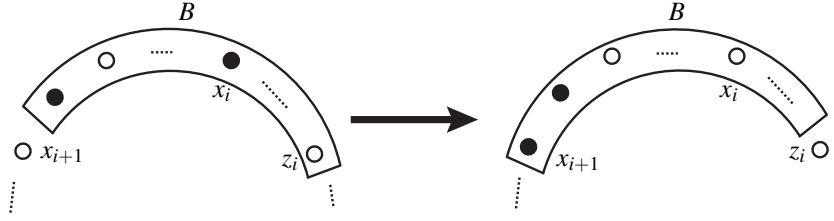
  Now note that we have the following property from our definition.
  \begin{align*}
    \size{X^{(i+1)} \cap D_0^{(i+1)}} &= \size{ \left( (D_0^{(i)}\backslash \{x_i\}) \cup \{x_{i+1}\} \right) \cap \left( (X^{(i)} \backslash \{x_i\}) \cup \{z_i\} \right) } \\
    &=\size{ \left( (D_0^{(i)} \cap X^{(i)}) \backslash \{x_i\} \right) \cup \left( D_0^{(i)} \cap \{z_i\} \right) \cup \left(X^{(i)} \cap \{x_{i+1}\} \right) } \\
    &= \size{D_0^{(i)} \cap X^{(i)} } - 1 + 0 + \size{X^{(i)} \cap \{x_{i+1}\} }
  \end{align*}
  For $0\leq i < s$, we know that $x_{i+1}\notin D_k$, and therefore
  $x_{i+1}\in X^{(i)}$. From the computation above, we see that this
  implies $|X^{(i+1)} \cap D_0^{(i+1)}| = |X^{(i)} \cap
  D_0^{(i)}|$. Therefore, the pair $(X^{(i+1)},D_0^{(i+1)})$ is
  optimal for $i<s$. On the other hand, $x_{s+1}\in D_k$, and
  therefore $x_{s+1}\notin X^{(i)}$, so we have $|X^{(s+1)} \cap
  D_0^{(s+1)}| < |X^{(0)} \cap D_0^{(0)}|$, contrary to the optimality
  of $(X^{(0)}, D_0^{(0)})$. Therefore, a set $D_k$ as described above
  cannot exist.
  \end{proof}

With this result in hand, we are now prepared to resolve the case
$c=4$.

\begin{theorem} \label{c=4} 
 Let $d$ be a positive integer, and let $n = 4d + 3$.  Then there
  exists a partition of the subsets of $[n]$ of size at least $d$ into
  intervals $[A,B]$ such that $|B| \geq d+3$.
\end{theorem}
\begin{proof}
Let $\mathcal{I} = \{[A, f_4(A)]\mid A \subseteq [n], |A| = d\}$.  By Lemma \ref{right size}, $|f_4(A)| = d+3$ for all $d$-sets $A$.  By Lemma \ref{d+1}, all $(d+1)$-sets are covered by the elements of $\mathcal{I}$.  Since Lemma \ref{right size} implies that $|f_4(A)|-|A| = 3$, we can apply Lemma \ref{no collisions} to the $d$-sets to show that $I \cap I' = \varnothing$ for all $I \not= I' \in \mathcal{I}$.  Therefore we have partitioned all of the $d$-sets and $(d+1)$-sets into the appropriate intervals.  What remains to show is that the $(d+2)$-sets and $(d+3)$-sets that are not covered by any element of $\mathcal{I}$ can be partitioned according to the theorem.

For $i \in \{1,2\}$ let $\mathcal{V}_i = \{S \subseteq [n]\mid |S| =
d+i+1, ~S \notin \cup_{[A,B] \in \mathcal{I}}[A,B] \}$.  Let $G$ be
the bipartite graph with vertex set $\mathcal{V}_1 \cup \mathcal{V}_2$
and $S_1S_2 \in E(G)$ if and only if $S_1 \in \mathcal{V}_1$ is a
subset of $S_2 \in \mathcal{V}_2$.  By Lemma \ref{uncovered top} with
$k=2$ we have that \[\Deg(S) = n-(d+2) \geq d+3\] for each $S \in
\mathcal{V}_1$.  Furthermore, $\Deg(S) \leq d+3$ for all $S \in
\mathcal{V}_2$.  Hence, by Corollary~\ref{cor:spec-hall}, there exists a
complete matching from $\mathcal{V}_1$ to $\mathcal{V}_2$.  Therefore
we can use the 1-cubes corresponding to the matching to cover the
elements of $\mathcal{V}_1$, and these intervals are disjoint.  The
remaining elements of $\mathcal{V}_2$ can then be covered by trivial
intervals.  Furthermore, the sets of size at least $d+3$ can be
covered by trivial intervals as well.  Thus, we have found a partition
of the desired type.  \end{proof}

This completes the proof of Theorem~\ref{thm:main}.

\section{Conclusion}

Although this paper only resolves a fraction of the possible cases of
Conjecture~\ref{CONJ-IND}, we believe that it is true in full
generality, in part because many of the lemmas in this paper are
stated and proved in greater generality than is required for their
applications here. However, it is not clear how to apply them to
advance this line of research. It seems that after applying
Lemma~\ref{uncovered top} it may be beneficial to pass from our
combinatorial setting back to the algebraic setting to acquire some
insight. Although it is not immediately obvious, our approach here
generalizes the constructive proof for $d=1$ given by Bir\'o et al.\ in
\cite{BHKTY}. The generalization comes through the notion of balanced
sets as originally defined there and an equivalent method of finding
the intervals. Unfortunately, we do not see how to generalize the
method of finding intervals that appears in \cite{BHKTY} to higher
values of $d$. We also note that the inductive proof for $d=1$ given
by Bir\'o et al.\ seems to face significant challenges in generalizing,
and while we would not rule out the use of such a method to prove
Conjecture~\ref{CONJ-IND}, it seems unlikely to be successful.

\section{Acknowledgments}

The authors are grateful to J\"urgen Herzog, David M.\ Howard and
William T.\ Trotter for useful conversations while pursuing this
research.


\end{document}